\documentclass[12pt]{amsart}
\usepackage[latin1]{inputenc}	% Accenti 

\usepackage[english]{babel}
\usepackage{amsmath,amsthm}
\usepackage{amssymb,latexsym}
\usepackage{graphicx}

\newcommand{\Z}{\mathbb{Z}}

\newtheorem{theorem}{Theorem}[section]
\newtheorem{lemma}{Lemma}[section]

\theoremstyle{remark}

\allowdisplaybreaks

\begin{document}
\title{Some remarks in elementary prime number theory}

\author{Raffaele Marcovecchio}
\address{Dipartimento di Ingegneria e Geologia\\
Universit\`a di Chieti-Pescara \\
Viale Pindaro, 42\\
65127 Pescara\\
Italy}
\email{raffaele.marcovecchio@unich.it}

\date{\today}

% AMS subject classifications (used in AMS journals)
   \subjclass[2010]{Primary ; Secondary }

% AMS keywords (used in AMS journals)
   \keywords{keywords}

\maketitle
\section{Introduction}
One of the many nice features of the Selberg integral (see ~\cite{Selberg}), i.e. of
\begin{equation}				\label{selberg}
\begin{split}				
S_n(\alpha,\beta;\gamma) & = 
	\int\limits_0^1\cdots\int\limits_0^1 
	\prod_{i=1}^n x_i^{\alpha-1} (1-x_i)^{\beta-1} 
	\prod_{1\leqslant i<j\leqslant n} |x_i-x_j|^{2\gamma} 
	{\rm d} x_1 \cdots {\rm d} x_n \\ & = 
	\prod_{j=0}^{n-1} \frac{\Gamma(\alpha+j\gamma) \Gamma(\beta+j\gamma)
	\Gamma(1+(j+1)\gamma)}{\Gamma(\alpha+\beta+(n+j-1)\gamma)\Gamma(1+\gamma)},
\end{split}
\end{equation}
is that it can be used in elementary methods for the study of the 
distribution of prime numbers. In the most general setting 
$\alpha$, $\beta$, $\gamma$ are complex numbers such that
\[
\mathrm{Re} (\alpha), \mathrm{Re} (\beta) > 0, \quad
\mathrm{Re} (\gamma) > - \min \left\{ \frac{1}{n}, 
		\frac{\mathrm{Re} (\alpha)}{n-1}, \frac{\mathrm{Re} (\beta)}{n-1}\right\}.
\]
Such elementary methods in prime number theory are developements 
of a method introduced by Gelfond and Schnirelmann~\cite{Gelfond} 
(see also~\cite[Chapter 10]{Montgomery} for a detailed analysis 
on this and some other related methods). We refer the reader 
to~\cite{ForresterWarnaar} for a survey on the relevance of the 
integral (\ref{selberg}), and to~\cite{Diamond} for a general 
survey on the prime number theory.

For $\gamma=1$ the integral (\ref{selberg}) evaluates a determinant of Hankel's type:
\begin{equation}			\label{determinantNair}
\frac{1}{n!} S_n(\alpha,\beta;1) = \det_{1\leqslant i,j\leqslant n} 
					\left( \int\limits_0^1 x^{\alpha+i+j-3} (1-x)^{\beta-1}
					{\rm d} x \right).
\end{equation}
This is easily seen by using a classical argument due to Heine 
(see e.g.~\cite{Zudilin}). In the papers~\cite{Nair} and~\cite{Chudnovsky}, 
independently, (\ref{determinantNair}) was used in order to generalize the 
Gelfond--Schnirelmann method. In~\cite{Nair} (and in~\cite{Chudnovsky} 
with similar considerations), by taking $\alpha=\beta=[sn]$ the greatest integer 
not exceeding $sn$, for $n\to\infty$ and $s=0.39191162...$, a new elementary proof 
of the following theorem is obtained:
\begin{theorem}
[Nair~\cite{Nair}]
For all sufficiently large $x$, we have
\begin{equation}			\label{chudnovskynair}
\psi_1(x) \geqslant (0.49517\dots) x^2 + O(x \log^2 x). 
\end{equation}
\end{theorem}
Here $\psi_1(x)$ is the sum function 
\[
\psi_1(x)= \sum\limits_{n\le x} \psi(n),
\]
and $\psi(x)$ is Chebyshev's $\psi$--function. Combining (\ref{selberg}) and
(\ref{determinantNair}), and evaluating the Euler beta integrals, we obtain
\begin{equation}				\label{mainformula}
\Gamma(\beta)^n \det_{1\leqslant i,j \leqslant n} 
		\left( \frac{\Gamma(\alpha+i+j-2)}{\Gamma(\alpha+\beta+i+j-2)} \right) = 
\prod_{j=0}^{n-1} \frac{\Gamma(\alpha+j) \Gamma(\beta+j)
							j!}{\Gamma(\alpha+\beta+n+j-1)}.
\end{equation}
In~\cite{Nair}, the required upper bound of $S_n([sn],[sn];1)$ is obtained 
by computing the maximum of the integrand in $S_n([sn],[sn];1)$ over $[0,1]^n$, 
while no use is made of the evaluation (\ref{mainformula}).

The aim of this paper is twofold. We give a somehow more direct proof of 
(a generalisation of) (\ref{mainformula}), and use this formula to obtain 
the required upper bound when $\alpha=\beta=[sn]$, thus getting a slightly 
new proof of (\ref{chudnovskynair}). Though none of these two remarks has special 
novelty, it seems hopeful that the general setting alluded to above, and even 
further generalizations occurring when $\gamma$ is a general positive integer 
(not necessarily $1$), might eventually lead to interesting developments. 
%%%%%%%%%%%%%%%%%%%%%%%%%%%%%%%%%%%%%%%%%%%%%%%%%%%%%%%%%%%%%%%%%%%%%%%
\section{A nice determinant}
%%%%%%%%%%%%%%%%%%%%%%%%%%%%%%%%%%%%%%%%%%%%%%%%%%%%%%%%%%%%%%%%%%%%%%%
In this section we review the proof of (\ref{chudnovskynair}) 
in~\cite{Chudnovsky} and~\cite{Nair}. We recall the following evaluation 
from~\cite{Krattenthaler}:
\begin{lemma} \cite[Lemma 3]{Krattenthaler}
For all indeterminates $X_1,\dots,X_n$, $A_2,\dots,A_n$, $B_2,\dots,B_n$ we have
\[
\det_{1\leqslant i,j\leqslant n}
	\Big( (X_i+B_2)\cdots(X_i+B_j)(X_i+A_{j+1})\cdots(X_i+A_n) \Big) 
		= \prod_{1\leqslant i<j\leqslant n} (X_i-X_j)
			\prod_{2\leqslant i\leqslant j\leqslant n} (B_i-A_j).
\]
\end{lemma}
By choosing $X_i=i$ ($i=1,\dots,n$), $B_i=\alpha+i-3$ and $A_i=\alpha+\beta+i-3$ 
($i=2,\dots,n$) one easily gets (\ref{mainformula}). 

From the Stirling formula it easily follows that
\[
\log (1!\cdots n!) = \frac{n^2}{2}\log n - \frac{3}{4} n^2 + O(n\log n) 
														\quad (n\to\infty).
\]
Let us denote by $\Delta_n(s)$ the quantity in (\ref{mainformula}) with 
$\alpha=\beta=[sn]$, where $s$ is a positive parameter to be chosen later. Then
\[
\Delta_n(s) = 1!\cdots (n-1)! \left( \frac{1!\cdots ([sn]+n-2)!}{1!\cdots ([sn]-2)!} 
								\right)^2 
				\frac{1!\cdots (2[sn]+n-3)!}{1!\cdots (2[sn]+2n-3)!}.
\]
Therefore
\begin{align*}
\log \Delta_n (s) & = \left( \frac{(2s+1)^2}{2}\log(2s+1)-s^2\log s 
		- (s+1)^2 \log(s+1) - 2(s+1)^2 \log 2 \right) n^2 			\\ 
				& + O(n\log n) \qquad (n\to\infty).
\end{align*} 
On the other hand, using the Pochhammer symbol $(x)_\beta$ for the shifted factorial, 
\[
\frac{\Gamma(\beta) \Gamma(\alpha+i+j-2)}{\Gamma(\alpha+\beta+i+j-2)} = 
\frac{(\beta-1)!}{(\alpha+i+j-2)_\beta} = 
\sum_{k=0}^{\beta-1} \binom{\beta-1}{k} \frac{(-1)^k}{\alpha+i+j+k-2}.
\]
It follows that
\[
d_{\alpha+\beta+i+j-1} \frac{(\beta-1)!}{(\alpha+i+j-2)_\beta}\in\Z, 
\]
where $d_m$ denotes the least common multiple of $1,\dots,m$. Hence 
\begin{equation}			\label{basicEstrapolation}
d_{\alpha+\beta+i+j-1} \frac{(\beta-1)!}{(\alpha+i+j-2)_\beta}\geqslant 1.
\end{equation}
However, due to the above determinant calculation,  we get the following 
improved inequality:
\begin{equation}			\label{improvedEstrapolation}
\prod_{i=1}^n d_{\alpha+\beta+n+i-3} 
	\frac{(n-i)! (\beta+i-2)!}{(\alpha+i-1)_{\beta+n-1}}	\geqslant 1.
\end{equation}
Therefore
\[
\psi_1(2[sn]+2n) - \psi_1(2[sn]+n) \geqslant - \log \Delta_n(s).
\]
Thus, after some calculations one arrives at (\ref{chudnovskynair}).

One may get in the mood of using this special case of \cite[Lemma 9]{Krattenthaler},
i.e.  \cite[Lemma 3]{Krattenthaler} with $X_i+B_j$ replaced 
by $(X_i+B_j)(X_i-B_j-C)$, $X_i+A_j$ by $(X_i+A_j)(X_i-A_j-C)$, $X_i-X_j$ 
by $(X_i-X_j)(C-X_i-X_j)$ and $B_i-A_i$ by $(B_i-A_j)(B_i+A_j+C)$. 
This can be easily obtained by putting $X_i=-(Y_i-C/2)^2$, $A_i=(E_i+C/2)^2$ 
and $B_i=(F_i+C/2)^2$.
%%%%%%%%%%%%%%%%%%%%%%%%%%%%%%%%%%%%%%%%%%%%%%%%%%%%%%%%%%%%%%%%%%%%%%%
\section{A pretentious generalization (but almost suitable for arithmetic progressions)}
%%%%%%%%%%%%%%%%%%%%%%%%%%%%%%%%%%%%%%%%%%%%%%%%%%%%%%%%%%%%%%%%%%%%%%%
If we let $B_i=i$ and $A_i=\beta+i$ in Lemma 2.1, we obtain the following 
generalization of (\ref{mainformula}):
\[
\det_{1\leqslant i,j\leqslant n} \left( \frac{(\beta-1)!}{(x_i+j+1)_\beta} \right) = 
	\frac{(\beta-1)!\cdots (\beta+n-2)!}{(x_1+2)_{\beta+n-1}\cdots (x_n+2)_{\beta+n-1}} 
	\prod_{1\leqslant i<j\leqslant n} (x_j-x_i).
\]
This implies that 
\[
d_{x_1+\beta+n}\cdots d_{x_n+\beta+n} 
	\frac{(\beta-1)!\cdots (\beta+n-2)!}{(x_1+2)_{\beta+n-1}\cdots (x_n+2)_{\beta+n-1}} 
	\prod_{1\leqslant i<j\leqslant n} (x_j-x_i) \geqslant 1,
\]
for all positive integers $\beta$, $x_1,\dots,x_n$ such that $x_1,\dots,x_n$ are distinct. 
The Chudnovsky--Nair determinant is the special case where $x_1,\dots,x_n$ are consecutive 
numbers.
%%%%%%%%%%%%%%%%%%%%%%%%%%%%%%%%%%%%%%%%%%%%%%%%%%%%%%%%%%%%%%%%%%%%%%%
\section{Conclusions}
%%%%%%%%%%%%%%%%%%%%%%%%%%%%%%%%%%%%%%%%%%%%%%%%%%%%%%%%%%%%%%%%%%%%%%%
No further interesting application seems to come out of considerations of 
hyperdeterminants \cite{LuqueThibon}, when $\gamma>1$, or $q$--analogs. 
Therefore the main contribution of this short note, if any, is to enlighten 
the improvement from (\ref{basicEstrapolation}) to (\ref{improvedEstrapolation}).
%%%%%%%%%%%%%%%%%%%%%%%%%%%%%%%%%%%%%%%%%%%%%%%%%%%%%%%%%%%%%%%%%%%%%%%

 %%%%%%%%%%%%%%%%%%
\end{document}